\documentclass[12pt,a4paper]{article}
\usepackage{amssymb,amsmath,amscd}
\usepackage[russian]{babel}
\usepackage[dvips]{graphicx}

\newtheorem{theorem}{\bf Theorem}
\newtheorem{Defin}{\bf Definition}
\newtheorem{remark}{\bf Remark}
\newtheorem{property}{\bf Property}

\usepackage{extsizes}

\begin{document}
\newcommand{\eqdefrl}{\stackrel{\mathrm{\Delta}}{=}}

\begin{center}
\textbf{\Large Skeleton decomposition of linear operators in the theory of degenerate differential equations\footnote{Joint work with Denis Sidorov and Yong Li. Partly supported by 
by the International science and technology cooperation program of China under Grant No. 2015DFR70850.}}\\
$\,$\\
{\it Nikolay Sidorov}
\end{center}

\textbf{Abstract. } {\it We suggest method based on the skeleton decomposition of linear operators in order to reduce ill-posed degenerate differential equation to the non-classic initial-value problem  enjoying  unique solution.}\\


\section{Problem statement}

This brief paper concerns Cauchy problem
\begin{equation*}
   \left\{ \begin{array}{lr}
         \mbox{$B \frac{dx}{dt} = x(t)+f(t)$} \hspace*{9cm} (1)\\
         \mbox{$x(0)=x_0, $} \hspace*{10.4cm} (2)\\
        \end{array} \right. 
\end{equation*}
where $B: X \rightarrow X$ is non-invertible linear operator, $X$ is banach space,
e.g. $B$ is Fredholm operator, $\ker B \ne \{0\}.$
The Cauchy problem  (1)-(2) is not solvable for arbitrary $ x_0$ making it  ill-posed problem in general settings.
Let us outline here that high-dimensional differential-algebraic-equations (DAEs) are the special case of such problem. DAEs are in the core of electromagnetic models of power systems.

If (1) is solvable then it's necessary 
$$\langle x(0)+f(0),\psi \rangle =0,  \eqno{(3)}$$ where $B^*\psi=0.$
This problem has been addressed by many authors. The approach described in 
 [1, 2] appears to be  productive in practical applications.
Nevertheless, due to condition (3), Cauchy problem (1)-(2) 
is not solvable in most of the cases and it is ill-posed problem
is sense of theory [3]. Here readers may also refer to ch. 5 of  textbook [4] .
In monograph [1] it's demonstrated that in sufficiently general settings, 
Cauchy problem (1)-(2)  is solvable in class of generalized functions.

In the next section we suggest to employ another initial condition for eq. (1) and propose the constructive approach which is easy of implement. We constract operator 
$M\subset {\mathcal L}(X \rightarrow X_p),$ where $X_p$ is linear normed space.
Th. 1 and Th. 2 make it possible to implement structure regularization
of  solution of eq.  (1) with  non-invertible  operator $B$
 using the
special initial condition selection
$$ M x(t)\bigr |_{t=0} =c_0, \eqno{(4)}$$
where $c_0 \in X_0$ instead of conventional Cauchy condition.
IVP (1), (4)  enjoy unique classic solution crossing the hyperplanes specified in proposed condition (4).
Solution of such regularized IVP (1), (4)   depends continuously on selection of 
$c_0 \in X_p.$

\section{Structure Regularization}

Let us introduce linear operators
$$A_{2i}:X_{i-1}\rightarrow X_i,\;\;\;A_{2i-1}:X_i\rightarrow X_{i-1},\;\;\;B^i\eqdefrl A_{2i}A_{2i-1}:X_i\rightarrow X_i,\;\;\;i=\overline{1,p},$$
where $X_i$ are linear normed spaces, $X_0\eqdefrl X$.

Let

\noindent{\bf 1)} $B=A_1A_2$, $\;\;A_{2i}A_{2i-1}=A_{2i+1}A_{2i+2}$, $\;\;i=\overline{1,p}$.

\begin{property}
$i$-th operator of the sequence  $\{B^i\}_{i=0}^p$ 
is constructed  by permutation of  skeleton decomposition of operators of {\it i--}1th operator.
Here $B^0 \eqdefrl B.$
\end{property}

Let us we impose the following condition

\noindent{\bf 2)} operators $B^0, B^1,..., B^{n-1}$ are noninvertible, 
and $B^p:\,X_p\rightarrow X_p$ is conti\-nuously invertible or non zero.

\begin{Defin}
Let the sequence $\{B^i\}_{i=0}^p$ meets conditions~1)~and~2). 
Then this sequence we call attached skeletal chain of finite length of the operator $B$. 
Moreover, of  $B^p$ is invertible then we can that chain as regular chain and we call it
degenerate chain is  $B^p$ is zero operator.
\end{Defin}

\begin{property}
Any finite non-invertible  operator generates skeletal chain of finite length [6, ch. II, \S 7].
\end{property}

\begin{remark} If $B$ is invertible, then its skeletal chain is reducible to operator $B$ itself.
\end{remark}

If conditions~1)~and~2) are satisfied then we can introduce functions $x_i(t)=\prod\limits_{j=1}^i A_{2j}x_0(t)$, where $x_0(t)$ is solution of eq. (1), $x_i(t)=A_{2i}x_{i-1}$, $i=\overline{1,p}$.
Let $\{ B_i\}_{i=0}^p$ be regular skeleton chain. Then function
 $x_p(t)$ 
 satisfies the regular Cauchy problem 
\begin{equation*}
   \left\{ \begin{array}{lr}
         \mbox{$B^p\frac{dx_p}{dt}=x_p+\prod_{j=1}^p A_{2j}f(t),$} \\
         \mbox{$x_p(0)=c_0. $}\\
        \end{array} \right. 
\end{equation*}
Once we have functions $x_p(t)$ determined, the rest of the functions $x_{p-1},\dots,x_0$ 
we construct following recursion
$$x_i=-\prod_{j=1}^i A_{2j}f(t)+A_{2i}\frac{dx_{i+1}}{dt},\; i=p-1,\dots,1,$$
$$x_0=-f(t)+A_1\frac{dx_1}{dt},$$
where $x_0$ is solution of initial problem (1), (4).
\begin{theorem} If operator $B$ has regular skeleton chain of  length $p$, $f(t)$ is 
$p$ times differentiable then eq. (1) with initial condition $$\prod\limits_{j=1}^pA_{2j}x(t)|_{t=0}=c_0,$$
where $\;c_0\in X_p$ enjoys unique classic solution $x_0(t,c_0)$. 
If spectrum of operator $B^p $   lies in the left half-plane, then solution of 
homogeneous IVP 
\begin{equation*}
   \left\{ \begin{array}{lr}
         \mbox{$B\frac{dx}{dt}=x,$} \\
         \mbox{$M x(t)|_{t=0}=c_0 $}\\
        \end{array} \right. 
\end{equation*}
 will be
asymptotically stable.  Here $M:=\prod\limits_{j=1}^pA_{2j},\, 0\leq t < \infty.$
\end{theorem}

\begin{theorem}

If operator $B$ has degenerate skeleton chain, i.e.
$B^p$ are zero operator, then homogeneous equation
$B\frac{dx}{dt} =x$ has only trivial solution and
unique solution of non-homogeneous equation (1)
can be constructed  as follows:
$$x(t) = -f(t)+A_1 \frac{dx_1}{dt},$$
where $x_1(t)$ can be constructed as following recursion
$$x_p(t) = -\prod_{j=1}^p A_{2j} f(t), $$
$$ x_i(t) = -\prod_{j=1}^i A_{2j} f(t) + A_{2i} \frac{dx_{i+1}}{dt}, \, i=p-1,\dots,1. $$
\end{theorem}

\begin{remark}
If operator $B$ has degenerate  skeleton chain
then $B$ is nilpotent operator and unique solution of nonhomogenious eq. (1)
can be also constructed as iteration $u_n(t) = -f(t)+B \frac{d}{dt} u_{n-1}(t), \, n=1,\dots, p,$
$u_0(t)=0,$ $u_p(t)$ satisfies a given eq. (1). 
\end{remark}

Such regularization scheme based on skeleton decomposition is applied in  [5] 
for the construction of trajectories passing through the hyperplanes $M x(0) = c_0$.\\

\noindent{\it Possible Generalizations.}\\
As footnote let us outline that similar results are valid for degenerate nonlinear Volterra
integral-differential equations
$$B \frac{d x(t)}{dt} = u(t)+\int_{0}^t K(t,s) g(s,x(s))\, ds,  $$
where $B$ is noninvertible linear operator. Of particular interest is the following case
$$ B = \sum_{i=1}^N a_i(t) \langle \cdot, \gamma_i \rangle, $$
where $\gamma_i \in X^*, \, a_i(t): {\mathbb R}^n \rightarrow X.$

$$\,$$

Results of this paper have been presented on the Intl Workshop on Inverse and Ill-posed Problems,
at Moscow State University (19-21 Nov. 2015) [7, p.165].

\begin{center}
\textbf{References}
\end{center}

\noindent 1. N.~Sidorov, B.~Loginov et al. Lyapunov-Schmidt methods in nonlinear analysis and applications. Springer Publ. 2013.

\noindent  2. G. Svridyuk, S. Zagrebina. Showalter-Sidorov problem as  phenomenon of Sobolev Equations. Bull. of Irkutsk State University, Ser. Maths, 2010. Vol.3, No. 1. pp.~104--125.

\noindent 3. A. N. Tikhonov, V. Y. Arsenin. Solutions of ill-posed problems. John Wiley \& Sons Inc. 1977.

\noindent 4. V. A. Tr\'enoguine. Analyse fonctionnelle. Moscou:  \'Editions Mir, 1985.

\noindent  5. O.A. Romanova, N.A. Sidorov. Trajectories construction of one dynamical system 
with initial data on hyperplanes. Bull. of Irkutsk State University, Ser. Maths, 2015. Vol.11. p.~93--105.

\noindent 6. F. R. Gantmacher. The theory of matrices. By F. R. Gantmacher. Trans. from  Russian by K. A. Hirsch, vols. I and II. New York, Chelsea, 1959.

\noindent 7. Abstracts of the Intl. Workshop on Inv. and Ill-posed Problems. 19-21 Nov. 2015. -- Moscow: Rus. Univ. Peoples Friendship Publ., 2015. - 217 p. ISBN 978-5-209-06803-7.

\end{document}